\documentclass[12pt]{article} 
\usepackage{amsmath,amsfonts,amssymb,amsthm}
\newcommand{\bb}{\mathbb}

\newcommand{\natls}{{\bb N}}

\newlength{\figboxwidth}             
\renewcommand{\bold}[1]{\medskip \noindent {\bf #1 }\nopagebreak}

\newcommand{\st}{\;\: : \;\:}         %Such that

\newcommand{\Tr}{\operatorname{Tr}}
\newcommand{\tr}{\operatorname{tr}}

\newcommand{\Ad}{\operatorname{Ad}}

\def\@ifundefined#1#2#3%
  {\expandafter\ifx\csname#1\endcsname\relax#2\else#3\fi}

\@ifundefined{theoremstyle}{
}{
\theoremstyle{plain} %default
}
\newtheorem{theorem}{Theorem}[section]

\newtheorem{lemma}[theorem]{Lemma}

\newtheorem{corollary}[theorem]{Corollary}

\@ifundefined{theoremstyle}{
}{
\theoremstyle{definition} %default
}
\newtheorem{definition}[theorem]{Definition}

\newcommand{\cC}{{\cal C}}

\newcommand{\cE}{{\cal E}}

\newcommand{\cO}{{\cal O}}

\newcommand{\cS}{{\cal S}}

\mathchardef\GG="321D
\title{Uniform Exponential Growth for Linear Groups}
\author{Alex Eskin\thanks{Research partially supported by the Packard
    Foundation.}, Shahar Mozes\thanks{Research partially supported by
    the BSF.} and Hee Oh\thanks{Research partially
    supported by NSF grant DMS0070544.} }

\begin{document}
\maketitle

\def\CC{{\mathbb C}}
\def\QQ{{\mathbb Q}}
\def\RR{{\mathbb R}}
\def\ZZ{{\mathbb Z}}

\def\d{{\rm d}}

\def\cC{{\cal C}}
\def\cO{{\cal O}}
\def\cS{{\cal S}}

\def\cOKS{{\cO_K(\cS)}}

\def\ueg{uniform exponential growth}

\section{Introduction}
Let $\Gamma$ be a finitely generated group. Given a finite set of
generators $S$ one has the corresponding  Cayley graph $\cC(\Gamma,S)$
and a word
metric $\d_S$ on $\Gamma$.
Let us denote by
$$B_S(n)= \{ \gamma\in \Gamma \st \d_S(\gamma,e)\le n\}, $$
the ball of radius $n$ in the Cayley graph $\cC(\Gamma,S)$.

\begin{definition}
\label{def:exponential:growth}
We shall say that $\Gamma$ has exponential growth if
$$\omega_S(\Gamma) \equiv \lim_{n \to \infty} |B_S(n)|^{1/n} > 1,$$ 
or equivalently,  $\lim_{n \to \infty} \tfrac{1}{n} \log |B_S(n)| >0$.
\end{definition}
Note that if $\Gamma$ has this property with respect to some generating set
$S$, it has it with respect to an arbitrary set of generators. 

\begin{definition}
\label{def:uniform:exponetial:growth}
We shall say that $\Gamma$ has uniform exponential growth if 
$$\inf_S \omega_S(\Gamma) >1,$$ 
where the infimum is taken over all finite
generating sets $S$ of $\Gamma$.
Equivalently $\Gamma$ has uniform exponential growth if and only if
there is some $c>1$ so that for any generating set $S$ and all $n \ge 0$ 
we have $|B_S(n)| \ge c^n$.
\end{definition}

A central question concerning groups of exponential growth which goes
back to the book \cite{Gr-La-Pa}
is whether a group $\Gamma$ of exponential growth is necessarily also
of uniform exponential growth. This was shown to be the case when the
group $\Gamma$ is hyperbolic by M.~Koubi \cite{Ko} (see also
\cite{De})  and in the case when $\Gamma$ is solvable
independently by D. Osin \cite{Os} and J. Wilson
\cite{Wi}. Recently R.~Alperin and G.~Noskov \cite{Alperin:Noskov:GL2}
have announced an
affirmative answer for 
certain subgroups of $SL_2(\CC)$. For a general discussion of these questions
see the survey \cite{DelaHarpe:Grigorchyk:survey}.  Our main result
is the following:
\begin{theorem}
Let $\Gamma$ be a finitely generated group which is linear over a
field of characteristic $0$. Then $\Gamma$ has exponential growth if
and only if it has uniform exponential growth.
\end{theorem}
Recalling that a linear group is of exponential 
growth if and only if it is not virtually nilpotent (see
\cite{Tits:pingpong}, \cite{Milnor}, \cite{Wolf})
we have an equivalent formulation:
\begin{corollary}
  If $\Gamma$ is a finitely generated group which is linear over a
  field of characteristic $0$, then either $\Gamma$ has uniform
  exponential growth or it is virtually nilpotent.
\end{corollary}

\section{Basic Observations}
Note that if a homomorphic image of a group has uniform exponential
growth then so does the original group.  Also if a finite index
subgroup has uniform exponential growth, then so does the original
group. Thus in view of the theorem
of Osin and Wilson mentioned above, we may assume that the Zariski closure
of $\Gamma$ is connected and simple, both in the algebraic sense.

\bold{Specialization.}  Let $E$ be the field of coefficients of $\Gamma$.
Note that since $\Gamma$ is finitely generated it follows that $E$ is
finitely generated.  Using the fact (see \cite{Lu-Ma}) that if a
finitely generated subgroup $\Lambda$ of $GL_n(\CC)$ is virtually
solvable then there is 
a bound on the index of a solvable subgroup in $\Lambda$ depending
only on $n$, we deduce that there exists a ``specialization'' i.e. a
field homomorphism $\sigma:E\to K$ inducing a homomorphism
$\rho:GL_n(E)\to GL_n(K)$ where $K$ is a finite extension of $\QQ$ so
that $\rho(\Gamma)$ is not virtually solvable. Hence we may assume
that we have a finitely generated group $\Gamma$ contained in
$SL_n(K)$ with $K$ a number field and having a simple connected
Zariski closure. 
Since $\Gamma$ is finitely generated, we may, after possibly replacing
$\Gamma$ by a 
finite index subgroup, assume that there is some finite set of
places $\cS$ of the field $K$ (which includes all the archimedean
ones) so that $\Gamma\subset SL_n(\cOKS)$, where $\cOKS$ denotes the
ring of $\cS$-integers in $K$. Thus the diagonal embedding
of $\Gamma$ in $\prod_{\nu\in\cS}SL_n(K_\nu)$ is discrete, where $K_\nu$
denotes the completion of $K$ with respect to the absolute value $|
\cdot |_\nu$ associated with the place $\nu \in \cS$. We choose,
for each $\nu$, an extension of the absolute value $|
\cdot |_\nu$ to the algebraic closure of $K$. 

\bold{Notation and terminology.}
Let us set the following convention of terminology which makes some of
the statements in the sequel more transparent. We shall use the term
``bounded'' to mean ``bounded with a bound depending only on the
dimension of the linear representation, the number field $K$ 
and the set of places $\cS$. 
We shall use the notation $x \prec y$ or equivalently $y \succ x$ to
mean that there exist bounded positive constants $c_1$ and $c_2$ such
that $x < c_1 y^{c_2}$. This is used only when $y \ge 1$. 
For a matrix $A \in SL_n(\bar{K})$, we define $\|A \|_\nu = 
\max_{ij} |A_{ij}|_\nu$, and $\|A \| = \max_{ \nu \in
  \cS} \|A\|_\nu$.  
\medskip

The following lemma plays a basic role and is one of the reasons why
we need to specialize so that the field of coefficients of $\Gamma$ is
a number field.
\begin{lemma}
\label{lemma:eigenvalues:separate}
For every $A \in SL_n(\cOKS)$, let $\lambda_1,\dots,\lambda_r$ be the
distinct eigenvalues of $A$. Then 
\begin{equation}
\label{eq:seperate}
\prod_{\nu\in \cS}\prod_{i\neq j}|\lambda_i-\lambda_j|_\nu \geq 1
\end{equation}
Hence for any $i\neq j$ and any $\nu\in\cS$ we have
$\frac{1}{|\lambda_i-\lambda_j|_\nu} \prec \| A \|$.
\end{lemma}

\bold{Proof.} The left hand side of (\ref{eq:seperate}) is the
absolute value of a non-zero (rational) integer, 
and is therefore at least $1$.
\qed\medskip 

In the course of the argument we shall need certain words of bounded
length with respect to any generating set of the group $\Gamma$ to be
as regular as possible. We shall use the following:
\begin{lemma}
\label{lemma:bezout}
Let $H$ be the connected Zariski closure of $\Gamma$ and let
$Z\subsetneq H$ be 
any Zariski closed proper subset of $H$. Then there exists $N=N(Z,H)
\ge 1$
such that for any finite generating set $S$ of $\Gamma$,
$B_S(N)\not\subset Z$. 
\end{lemma}

Lemma~\ref{lemma:bezout} is proved using the generalized Bezout
theorem.
\medskip

To show that our group $\Gamma$ has uniform exponential growth we
shall show that there is some bounded constant $m$  so that
given any finite generating set $S$ there exists a pair of elements in the
ball $B_S(m)$ generating a free semigroup. We recall the well-known
result of J.~Tits  \cite{Tits:pingpong} 
which states that any non-virtually-solvable linear
group contains two elements $A$ and $B$
which generate a free subgroup; the proof is based on the so called
``ping-pong lemma''. Our result
may be viewed as a sort of quantitative version of Tits' theorem, in
the sense that we obtain a uniform bound on the word length of the elements
$A$ and $B$; however our elements are only guaranteed to generate 
a free semigroup. 

Showing that a pair of elements generates a free semigroup is based on
the following version of the ping-pong lemma which is due to G.~A.~Margulis.

\begin{definition}[Ping-Pong Pair]
\label{def:ping-pong:pair}
Let $k$ be a local field. 
A pair of matrices $A,B\in SL_n(k)$ is a {\em ping-pong pair} if there
exists a nonempty set $U\subset{\mathbb P}(k^n)$ such that:
$$BU\cap U =\emptyset$$
$$ABU\subset U, \qquad A^2BU\subset U$$
\end{definition}

\begin{lemma}[Margulis]
\label{lemma:ping:pong}
If $A$, $B$ are a ping-pong pair then $AB$ and $A^2B$ generate a free
semigroup.
\end{lemma}

To apply this lemma we need an effective way of showing that certain
elements form a ping-pong pair.

\begin{lemma}
\label{lemma:inequalities:ping:pong}
Let $k$ be a local field with an absolute value $|\cdot |_\nu$. 
Suppose that $A,B\in SL_n(k)$ are matrices such that after
conjugation of both by a common matrix we have:
\begin{itemize}
\item[{\rm (L1)}]
$A=\begin{pmatrix}a_1&&&&\\
            &a_2&&&\\
            &&\ddots&&\\
            &&&&a_n
\end{pmatrix}
$ with $|a_1|_\nu\ge|a_2|_\nu\ge\cdot\ge|a_n|_\nu$ and
$|a_1|_\nu\ge\max(2,2|a_2|_\nu)$.
\item[{\rm (L2)}] There exist constants $c_2 >0$, $d_2 > 0$ such that
$\frac{1}{|B_{11}|_\nu} \le c_2 \|A\|_\nu^{d_2}$,  
and $B e_1 \not\in k e_1$.
\item[{\rm (L3)}] There exist constants $c_3 >0$, $d_3 > 0$ such that
$\|B\|_\nu \le c_3 \|A\|_\nu^{d_3}$.
\end{itemize}
Then there exists a constant $m$, depending only on $n, c_2, c_3, d_2$
and $d_3$,  such that $A^m$ and $B$ form a
ping-pong pair. 
\end{lemma}

\section{The steps of the proof}
\bold{Step 1.}
Using Lemma~\ref{lemma:bezout} one deduces that given any finite generating
set $S$ of $\Gamma$ there is a pair of elements $A, B$ of $\Gamma$ belonging to
a ball of bounded diameter in the Cayley graph $\cC(\Gamma,S)$ so that
any word of a certain very large but bounded length in $A$ and $B$ yields an
$H$-regular element, where $H$ denotes the Zariski closure of
$\Gamma$. (By an $H$-regular element we mean a semisimple element the
dimension of whose centralizer is equal to the rank of $H$. Such
elements exist in $H$ by \cite{Pr}.)
Moreover (using Lemma~\ref{lemma:bezout} for
$\Gamma\times\Gamma< H\times H$) we may require also that no two of
these bounded length words lie in a common parabolic subgroup of $H$.
We shall actually require these regularity and genericity conditions
also for all the corresponding words in $\rho_m(A)$ and $\rho_m(B)$
where $\rho_m$ denotes the $m$'th wedge representation of $SL_n$,
$1\le m \le n$. 
\medskip

For the next two steps, it is convenient to use the diagonal embedding
and consider $A$ and $B$ as
elements of $G = \prod_{\nu \in \cS} SL_n(\bar K_\nu)$. 
Also for a matrix $C = \prod_{\nu \in \cS} C_\nu \in G$, we write
$\| C \|_\nu = \max_{ij} | (C_\nu)_{ij} |_\nu$, and $\| C \| = \max_{
  \nu \in \cS } \| C \|_\nu$. 

\bold{Step 2.}
As $A$ is semisimple we may conjugate both $A$ and $B$ by the same matrix
in $G$ so that $A$ becomes diagonal. Conjugating by an appropriate element
from the centralizer of $A$ and possibly replacing $B$ by some word of
bounded length in $A$ and $B$ we may ensure that either:
%\begin{itemize}
%\item[{\rm (a)}]
\begin{equation}
\label{eq:B:prec:A}
\|B\| \prec \|A\|,
\end{equation}
%\end{itemize}
or,
%\begin{itemize}
%\item[{\rm (b)}] 
for some bounded $m$, 
\begin{equation}
\label{eq:trace:big}
\|A\|^m \max_{\nu \in \cS}|\Tr B|_\nu  \succ\|B\|.
\end{equation}
%\end{itemize}

\bold{Step 3.}
We claim we can modify our choice of $A$ and $B$ and conjugate by a common
matrix in $G$ so that (\ref{eq:B:prec:A}) holds and $A$ is diagonal.
\medskip

Indeed, suppose (\ref{eq:B:prec:A}) does not hold. Then, in view of
Step~2, (\ref{eq:trace:big}) holds. 
Also since (\ref{eq:B:prec:A}) does not hold, we may assume 
\begin{equation}
\label{eq:A:prec:B}
\|A \|^m \le \|B \|^{1/2},
\end{equation}
where $m$ is as in (\ref{eq:trace:big}). (Recall that $A$ is diagonal
and has determinant $1$, hence $\|A \| \ge 1$.)
Hence,  in view of 
Lemma~\ref{lemma:eigenvalues:separate} and the 
semisimplicity of $B$, 
we can diagonalize $B \in G$ using a matrix 
$C \in G$ such that
$\max\{\|C \|,\|C^{-1}\|\} \prec \|B\|$. 
Hence in view of (\ref{eq:A:prec:B}), 
\begin{equation}
\label{eq:tmp1}
\| C A C^{-1} \| \le \| C \| \|A \| \|C^{-1}\| \prec \| B\|. 
\end{equation}
Since $C B C^{-1}$ is diagonal, 
and in view of (\ref{eq:trace:big}) and (\ref{eq:A:prec:B}), 
\begin{equation}
\label{eq:tmp2}
\| C B C^{-1} \| \succ \max_{ \nu \in \cS} |\Tr B|_\nu \succ
\|A\|^{-m} \|B \| \ge \|B\|^{1/2}
\end{equation}
Hence in view of (\ref{eq:tmp1}) and (\ref{eq:tmp2}), 
$\|C B C^{-1}\| \succ \|C A C^{-1}\|$. Hence if (\ref{eq:B:prec:A})
fails  we can conjugate $A$ and $B$ by $C$  
and interchange the roles of $A$ and $B$.

\bold{Step 4.}
Let $\nu \in \cS$ be such that $\|A\|_\nu$ is maximal. 
Let $\rho_m$ be as in Step 1. From the discreteness of $SL_n(\cOKS)$, 
there exists $m$, $1 \le m \le n/2$ such that
$\rho_m(A_\nu)$ satisfies (L1).
Note that in view of (\ref{eq:B:prec:A})  we have (L3) for the elements 
$\rho_m(A_\nu)$ and $\rho_m(B_\nu)$.
Next we claim that by replacing $\rho_m(B_\nu)$ by a bounded word in
$\rho_m(A_\nu)$ and $\rho_m(B_\nu)$ we can ensure that (L2) is satisfied. 
Thus we produced a ping-pong pair using words of bounded length in the
given set of generators.

\section{Construction of words with nice properties}
Observe that in  steps 2 and 4 above we need to produce words of
bounded length so that the coefficients of the corresponding matrices
satisfy certain properties. Let us describe first how this is achieved
in the case where the Zariski closure of our
group $\Gamma$ is $SL_n$. The reason for this case being simpler than the
general one is that in this case we may assume that the diagonalized
matrix $\rho_m(A)$ has distinct eigenvalues. For the following, we fix
$\nu \in \cS$. Also let us for
simplicity of notation speak about $A$ and $B$ rather than about
$\rho_m(A)$ and $\rho_m(B)$, and write $| \cdot |$ and $\| \cdot \|$
rather than $| \cdot |_\nu$ and $\| \cdot \|_\nu$. 
For a constant $c>0$ define an
entry $B_{ij}$ of $B$ to be {\em large} when $|B_{ij}|> c \|B\|$.
Let ${\cal E}(c)$ denote the algebra generated by all the elementary
matrices $E_{ij}$ corresponding to large entries of $B$.  The
following lemma is used twice, once in the proof of Step 2, with $c$
being a small absolute constant, and once in the proof of Step 4,
with $c = \|A\|^{-m}$, with $m$ a large but bounded positive constant. 
\begin{lemma}
\label{lemma:vandermonde:large}
If $E_{ij} \in \cE(c)$ then there exists a bounded word  $C$ in $A$
and $B$ so that $|C_{ij}| > p \|A\|^{-k} \|B\|^l$ where $p$, $k$ and $l$ are
bounded positive constants.
\end{lemma}
Rather than giving here the proof of the lemma let us give an example
explaining the main idea. Assume that $B_{12}$ and $B_{21}$ are big
entries. Then we claim that some word of the form $BA^kB$ with $0\le k
\le n-1$ has big $(1,1)$ entry. Indeed we have 
\begin{equation}
\label{eq:matrix:mult}
(BA^kB)_{11}= \sum_i
B_{1i}B_{i1} a_i^k. 
\end{equation}
Note that $B_{12} B_{21}$ is a large term, but a
priori, it may cancel with other terms. However, it may not happen for
all $0 \le k \le n-1$, in view of Lemma~\ref{lemma:eigenvalues:separate}. 
More formally, we may argue as follows: 
Written vectorially, (\ref{eq:matrix:mult}) is
\begin{displaymath}
\begin{pmatrix}(BA^0B)_{11} \\ (BA^1B)_{11}
  \\ \vdots \\(BA^{n-1}B)_{11})\end{pmatrix} = \begin{pmatrix}
1 & a_1 & a_1^2 & \dots & a_1^{n-1} \\
1 & a_2 & a_2^2 & \dots & a_2^{n-1} \\
 \vdots \\
1 & a_n & a_n^2 & \dots & a_n^{n-1} \\
\end{pmatrix}
\begin{pmatrix} B_{11}^2 \\ B_{12}B_{21} \\ \vdots \\ B_{1n}B_{n1}
\end{pmatrix}
\end{displaymath}
Hence using the formula for the Vandermonde determinant and
Lemma~\ref{lemma:eigenvalues:separate} we deduce that since we have
$B_{12}B_{21}$ big we cannot have $(BA^kB)_{11}$ small for all $0\le k
\le n-1$.

\bold{Almost Algebras.}
In the general case, another technical device is needed. 
For matrices $A$ and $B$, 
let $\langle A, B \rangle  = \tr A B^{t}$, and let 
$\epsilon>0$ be a small parameter. We say that matrices $A_1, \dots
A_k$ form an {\em $\epsilon$-almost-algebra} if for each $1 \le i,j \le
k$, $\langle A_i, A_j \rangle = \delta_{ij}$, and $A_i A_j$ is
within $\epsilon$ of a linear combination of  $A_1, \dots, A_k$.  

For each $1 \le m \le n^2+2$, let us choose constants $\epsilon_m$ so that 
$\epsilon_{m+1} \ll \epsilon_m$, and $\epsilon_1 \ll \epsilon$. 
Now let $B_1, \dots B_m$ be matrices of norm at most $1$. 
For each $1 \le k \le n^2+1$, let 
\begin{multline*}
f(k) = \inf \{ m \in \natls | \text{ there exists a subspace $V_m$ of
  dimension $m$ such that all words } \\
  \text{ in the $B_i$ 
  of length at most
  $2^k$ are within $\epsilon_k$ of $V_m$. } \}
\end{multline*}
By construction, the function $f$ is increasing, and is bounded by
$n^2$. Hence there exists a minimal $k$, $1 \le k \le n^2+1$ such
that $f(k) = f(k+1)$. Then an orthonormal basis for $V_{f(k)}$ is an 
$\epsilon$-almost-algebra. We call it the 
{\em almost-algebra generated by $B_1, \dots, B_m$}.

In the proof we use, as a replacement to $\cE(c)$, 
the almost-algebra generated by the ``blocks'' of
$B$, (i.e. the projections of $B$ onto the eigenspaces of $\Ad(A)$).
It is easy to see that the analog of
Lemma~\ref{lemma:vandermonde:large} is satisfied. Also, 
in the appropriate sense, every almost algebra
is near an algebra. This allows us to complete the proof.

\vspace{1in}

\bigskip
Department of Mathematics, University of Chicago, Chicago IL~60637,
USA; \\
eskin@math.uchicago.edu

\bigskip
Institute of Mathematics, Hebrew University, Jerusalem~91904, ISRAEL;\\
mozes@math.huji.ac.il

\bigskip
Department of Mathematics, Princeton University, Princeton NJ, USA;\\
heeoh@math.princeton.edu

\end{document}